\documentclass[11pt]{amsart}\usepackage{epsfig}

\usepackage{amsmath}
\usepackage{amssymb}
\usepackage{amsthm}
\usepackage{epsfig}
\usepackage{graphicx}

\newcommand{\RR}{\rm I\kern -1.6pt{\rm R}}

\newtheorem{theorem}{Theorem}[section]
\newtheorem{lemma}{Lemma}[section]

\numberwithin{figure}{section} \numberwithin{equation}{section}
\begin{document}
\date{}
\title[A potential system with infinitely many critical periods]
{A potential system with infinitely many critical periods}

\author{\sc Jihua Wang}
\address{}
\email{wangjh2000@163.com (J. Wang)}
%\thanks{Research was supported by the National Natural Science Foundations of China (No.XXX)}

\keywords{~Period function; Potential system; Critical point;
Calculus of Variation.}

% It is required to enter 2010 MSC.
\subjclass{34C23, 34C25, 37G15}

\begin{abstract}
In this paper, we propose an analytical non-polynomial potential
system which has infinitely many critical periodic orbits in phase plane. The main approach is Rolle' Theorem by showing the existence of infinitely many $2\pi-$ periodic solutions, through utilizing variational methods and the properties of Bessel function. The result provides an affirmative example to Dumortier's conjecture [Nonlinear Anal.
\textbf{20} (1993)].
\end{abstract}
\maketitle

\section{Introduction and main result}
The phenomenon of Oscillations (or vibrations) is popular in
dynamical system. The study on its period (or the frequency) is very
important in the applications of engineer or the other technology.
Consider the conservative Newtonian model $\ddot{x}+g(x)=0$ or equivalently Hamiltonian system in the phase plane
\begin{equation}\label{SN}
\begin{split}
\frac{dx}{dt}=-y,\quad \frac{dy}{dt}=g(x).
\end{split}
\end{equation}
It generates the classical one degree of freedom Hamiltonian vector
fields with a ``kinetic + potential" Hamiltonian of the
form $H(x,y)=\frac{y^2}{2}+G(x)=h$, where $G(x)=\int_0^xg(s)ds$
stands for the potential energy.

Assume that the origin $O$ is the center of system \eqref{SN}, the
biggest punctured neighborhood enclosing which is called a period
annulus. Denote it by $\mathcal{P}$, is foliated by periodic orbits.
The projection of which onto $x-$axis by $(x_l,x_r)$.
Assigning to each periodic oval in $\mathcal{P}$ its minimal period of the motion
along it, then we get period function $T$ of the center. Since
the period function is defined on the set of periodic orbits in
$\mathcal{P}$, it is natural firstly to parameterize the set.
Usually parameterizing it by the energy level $h$ of the periodic
oval and denote the oval by $\Gamma_h$ accordingly. It is known that
for a given continuum of periodic orbits, the number, character
(maximum or minimum) and distribution of critical periods do not
depend on the particular parametrization of the set of periodic
orbits applied. If the period function $T(h)$ is constant in a
neighborhood of the center, the center is called isochronous. The
isolated zeros $\hat{h}$ of the first order derivative of $T$ are
called critical periods, accordingly the periodic oval
$\Gamma_{\hat{h}}$ is called critical periodic orbit. The study on
period function is closely related to two-point nonlinear boundary-value problem of ODEs \cite{CC}
and some bifurcation problems
such as subharmonic bifurcations of periodic oscillations and the
bifurcation of steady-state solutions of a reaction-diffusion
equation \cite{ZO,URABE,SCJS}.

As commented in \cite{JF,EFAA}, many problems related to the period
function and its critical point are the counterpart of the poincar\'{e} return map
and limit cycle (isolated periodic orbit), the latter is the main focus of the celebrated Hilbert's 16th
problem. In 1980, Z. Zhang \cite{ZZ} found an example that a planar
analytic system has an infinite number of limit cycles which
answered a conjecture opened for many years. where it was shown analytical property of
vector field only assure the local finiteness of number of limit
cycles for a planar system of differential equation while it could
not deduce the global finiteness of limit cycles.

It is analogous to non-accumulating of limit cycles in a compact set of analytic vector field. Chicone and Dumortier \cite{CCFD} proved that the critical periodic
orbits of an analytic center can not accumulate on a compact set. If
potential $g$ of System \eqref{SN} is nonlinear
polynomial, then the number of critical periods is finite globally.
To the best of our knowledge, there are few papers on the existence
of multiple (no less than two) critical periods
\cite{CCMJ,EFAA,AJCLJY,MGJV} which should be taken into
account when looking for multiple solutions of boundary value
problems \cite{MSTI}. F. Dumortier remarked in \cite{CCFD} ``it is possible
that an analytic vector field has an infinite number of critical
periods. For example, a linear vector field on the plane with a
period annulus has a constant period function". Nevertheless, limit cycle is commonly defined as the isolated periodic orbit, in present paper we intend to propose an nontrivial example
where a potential system has infinitely many critical periods demonstrated as the
strict extremum of period function.

Consider system \eqref{SN} with potential $g(x)=x+\lambda\sin(x)$.
\begin{equation}\label{SP}
\begin{split}
\frac{dx}{dt}=-y, \quad \frac{dy}{dt}=x+\lambda\sin(x),
\end{split}
\end{equation}
where the parameter $\lambda\in\,[-1,4.603338849)$ is real and the potential
energy
\[G(x)=\int_0^x\,g(s)ds=\frac{x^2}{2}+\lambda(1-\cos(x))\]
is even. For all $\lambda\in\,[-1,+1]$ the period annulus is the whole
plane, in this setting the origin is called a global center. The periodic orbit is
the level curve
$\{(x,y)\in\,R^2|\,\frac{y^2}{2}+G(x)=h,\;h\in\,(0,+\infty)\}.$ When
$\lambda\in(-1,+1],$ the center is elementary, the Jacobian Matrix at which has a pair of pure
imaginary eigenvalues, i.e., the linearized matrix is non-degenerate.
When $\lambda=-1$, the origin is a nilpotent center.

The system \eqref{SN} with potential $g(x)=x+\lambda\sin(x)$ could be regarded as
a perturbation of the rigidly isochronous center at the origin
(i.e., $\dot{\theta}\equiv\,1$ in the polar coordinates, as defined
in \cite{JCMS}). In another hand C. Chicone \cite{CC} consider the Neumann boundary value problem
\begin{equation}
\begin{split}
\ddot{x}+\mu\,\sin x=0,\\
\dot{x}(0)=0,\quad \dot{x}(T)=0.
\end{split}
\end{equation}
with $\mu>0$. Here $x\equiv\,0$ is always a solution. One may be interested in subharmonic solutions with $N$ nodes, that is $x(t)=0$ exactly $N$ times on $(0,T)$. In geometrical viewpoint the period annulus enclosing the origin has the minimum period $\frac{2\pi}{N}$. It reveals that the analysis of boundary value problem can be completely understood in terms of period function which assigns to each $x\in\,(0,\pi)$. It is closely related to the problem of critical periods.

In present paper, we are interested in the
bifurcation of critical periods, roughly speaking, the emergence or
disappearance of critical period by perturbing a period annulus
enclosing a center. In particular perturbing the periodic orbits surrounding isochronous center
\cite{EFAA,MGJV}, we ask for the number of critical period orbits that persists.
In a similar manner as in studying
limit cycle bifurcation whether if the period orbit persists after a
small perturbation of the continua of periodic orbits enclosing a
center, which is known as the weakened Hilbert's 16th problem
posed by Arnold \cite{VIAD}.

The main result is stated as follows
\begin{theorem}\label{MRT}
The origin of potential system \eqref{SP} has infinitely many critical
periods, it appears alternatively as maximum and minimum.
\end{theorem}

\section{Preliminary lemmas and results}
The asymptotic expansion of period function near the simple center
is well known, See Theorem 4.3 of \cite{CCFD} or Proposition 3.5 of
\cite{LYXZ} for instance.
\begin{lemma}\label{PHDPH}
If $\lambda\in\,[-1,\lambda_m)$, the origin $O(0,0)$ of potential system
\eqref{SP} is a global center, where $\lambda_m=4.603338849)$.

When $\lambda\in\,(-1,\lambda_m)$, the center is simple and its period function $T(h)$ can be
analytically extended to the origin, the period function and its
first order derivation $T^{\prime}(h)$
\begin{equation}
         T(0^+)=\frac{2\pi}{\sqrt{g^{\prime}(0)}},\;
T^{\prime}(0^+)=\pi\frac{{5(g^{\prime\prime}(x))^2-3\,g^{\prime}(x)g^{\prime\prime\prime}(x)}}{12(g^{\prime}(x))^{\frac{7}{2}}}|_{x=0}.
\end{equation}

When $\lambda=-1,$ the origin $O(0,0)$ is the nilpotent center of
order 1, the period function $T$ has a pole of order one at $O$, that
is $T(0^+)=+\infty$ and $T^{\prime}(0^+)=-\infty.$
\end{lemma}

\begin{proof}
By linearization, it is easy to see the origin of planar Hamiltonian system \eqref{SP} is simple center as $\lambda>-1$ and nilpotent center as $\lambda=-1$, respectively.
Furthermore we consider the conditions which ensure the origin $O(0,0)$ of potential system
\eqref{SP} to be a global center, that is potential function $g(x)=x+\lambda\,\sin(x)$ has no positive zeros since it is odd. It holds $|\sin(x)|\leq\,|x|$ for all real $x$, thus for all $\lambda\in\,[-1,1]$, potential $g(x)$ has a unique zero at $x=0$. Therefore we aim to determine the minimal parameter $\lambda>1$ under which $g(x)$ has positive zeros. Noticing $g(x)>0$ for all $x\in\,(0,\pi)$ and $\frac{\partial\,g(x)}{\partial\,\lambda}=\sin(x)<0$ as $x\in\,(\pi,2\pi)$. The minimal positive parameter $\lambda_m>1$ suffices the conditions
\[g(x)=x+\lambda\,\sin(x)=0,\quad g^{\prime}(x)=1+\lambda\,\cos(x)=0\]
as $x\in\,(\pi,2\pi).$ In geometrical view of point, it implies that the smooth curve $g(x)$ tangentially intersects the positive half $x-$axis. It deduces that transcendental equation $\tan(x)=x$ has a unique solution $x\approx\,4.493409458$ at the interval $(\pi,2\pi)$. From which we obtain an approximate value $\lambda_m\approx\,= 4.603338849.$
\end{proof}

In the setting of global center, concerning the asymptotical behaviors of $T(h)$ near
the outer boundary of $\mathcal{P}$, we have
\begin{lemma}(\cite{DTHHZF})\label{LEMOB}
If potential $g$ of system \eqref{SN} is locally Lipschitz and
satisfies
\begin{enumerate}
  \item [(1).] strict condition of signs:~$xg(x)>0$ for all $x\neq\,0$;
  \item [(2).] semilinear condition:~$0<\liminf\limits_{|x|\rightarrow\,+\infty}\frac{g(x)}{x}\leq\,\limsup\limits_{|x|\rightarrow\,+\infty}\frac{g(x)}{x}<+\infty,$
\end{enumerate}
then the period function $T(h)$ satisfies
\[0<\liminf\limits_{h\rightarrow\,+\infty}T(h)\leq\,\limsup\limits_{h\rightarrow\,+\infty}T(h)<+\infty.\]
\end{lemma}

When the potential $g$ of planar system \eqref{SN} is odd, Concerning the analytical properties of period function of potential center. Urabe's striking conclusion \cite{URABE} concerning its iso-chronicity states as follows
\begin{lemma}\label{URABE}
If the potential $g$ is odd function, the only type of planar 
potential system \eqref{SN} having an isochronous center at the
origin is the linear one.
\end{lemma}

Concerning the monotonicity of period for potential center \eqref{SN} at the origin, Opial assert in Theorem 6 in \cite{ZO}
\begin{lemma}\label{OLEMA}
If the potential $g$ is odd function, the period function of global center of system \eqref{SN} is
monotonous then the function
$x\mapsto\,\frac{G(x)}{x^2}$ is monotonous in $(0,+\infty)$ as well.
\end{lemma}

Parameterizing the periodic orbit by the abscissa $\xi$ of
intersection point between oval $\Gamma$ with positive half
$x-$axis, denoted by $\Gamma_{\xi}$. Variation Lemma 3.2 of \cite{CCMJ} stated as follows proposes a necessary condition to the persistence of $2\pi-$periodic solutions.
\begin{lemma}\label{VLEMA}
Let $(x,y)\rightarrow\,(f(x,y),g(x,y))$ be analytic mapping,
consider the one parameter family of perturbed system from the
linear center as follows
\[\frac{dx}{dt}=-y+\lambda\,f(x,y),\quad \frac{dy}{dt}=x+\lambda\,g(x,y),\]
which has a center at the origin for any parameter $\lambda\in\,J$ where real closed interval $J$ contains $\lambda=0$.
Then period function $T(\xi,\lambda)$ holds for
$T_{\lambda}^{\prime}(\xi,0)=0$ if and only if the integral
$$\int_0^{2\pi}[xg(x,y)-yf(x,y)]_{x=\xi\cos\theta,y=\xi\sin\theta}\,d\theta=0.$$
\end{lemma}

{\bf proof.}~The conclusion (i) is proved in \cite{CCMJ} by using variation of parameter. Note that the origin $O$ is an linear isochronous center as $\lambda=0$ and the period function $T(\xi,0)\equiv\,2\pi.$
By the continuous
dependence on the initial value and parameter associated with the
solution of differential equation, the period of closed orbit
$\Gamma_{\xi,\lambda}$ is the variation of $2\pi$, denoted by
$T(\xi,\lambda)=2\pi+\rho(\xi,\lambda)$. It is clear $\rho(\xi,0)=0$ and $\frac{\partial{\rho(\xi,0)}}{\partial{\lambda}}=0$ is equivalent to vanishing of the following integral
$$\int_0^{2\pi}[xg(x,y)-yf(x,y)]_{x=\xi\cos\theta,y=\xi\sin\theta}\,d\theta=0.$$
That ends the proof.

The contour integral of perturbation terms $\lambda[f(x,y)\frac{\partial}{\partial\,x}+g(x,y)\frac{\partial}{\partial\,y}]$ along the oval passing through $(\xi,0)$ vanishes which proposes a necessary condition to persistence of $2\pi$ periodic solutions after perturbations.

\section{The proof of Theorem \ref{MRT}}
In this section we intend to prove Theorem \ref{MRT} by twofold: the
first one is the existence of critical periods, another one is the
unboundedness of number of critical periods.

Firstly it is obvious when $\lambda\neq\,0$ the origin of system \eqref{SP}
is non-isochronous center as shown in Lemma \ref{URABE}.

Using Lemma \ref{PHDPH} we get the asymptotic expansion of
$T(h,\lambda)$ for $h$ sufficiently small is
\begin{equation}\label{ASPTH}
T(h,\lambda)=\frac{2\pi}{\sqrt{1+\lambda}}+\frac{\lambda\,\pi}{4(1+\lambda)^{\frac{5}{2}}}\,h+O(h^2).
\end{equation}

It follows from Lemma \ref{LEMOB} that the system \eqref{SP} is
semilinear, hence period function $T(h)$ for $h\rightarrow\,\infty$
satisfies
\[0<\liminf\limits_{h\rightarrow\,+\infty}T(h)\leq\,\limsup\limits_{h\rightarrow\,+\infty}T(h)<+\infty.\]

By direct computation we get for system
\eqref{SP} that
\[\frac{G(x)}{x^2}=\frac{1}{2}+\frac{\lambda(1-\cos(x))}{x^2},\]
its derivative with respect to $x$ is
\[[\frac{G(x)}{x^2}]^{\prime}(x)=\lambda\frac{x\sin(x)-2+2\cos(x)}{x^3}.\]
An easy computation shows that function $\frac{G(x)}{x^2}$ is
monotonous in $(0,2\pi)$, but not monotonous in $(0,\infty)$, it follows from Lemma \ref{OLEMA} the period function of system \eqref{SN} is not monotonous
in $(0,+\infty)$. Therefore potential system \eqref{SP} has critical periodic orbits
in the unbounded period annulus $R^2$.

By Lemma \ref{VLEMA}, we intend to consider the variational problem on the period
$T(\xi,\lambda)$. By change of polar coordinates, we get for the $2\pi-$periodic orbit
$\Gamma_\xi$ of system \eqref{SP} it satisfies
\begin{eqnarray*}
\int_0^{2\pi}\sin(\xi\cos(s))\cos(s)\,ds&=&4\int_{0}^{\frac{\pi}{2}}\sin(\xi\sin(u))\sin(u)\,du\\
                                         &=&4\int_0^1\sin(\xi\,w)\frac{w}{\sqrt{1-w^2}}\,dw
                                         =2\pi\cdot\,J(1,\xi)=0.
\end{eqnarray*}

Where the notation $J(1,\xi)$ stands for the Bessel function of the
first kind \cite{GRR} which satisfies the following Bessel's
equation of order one
\[\xi^2\frac{d^2\,x}{d\xi^2}+\xi\frac{dx}{d\xi}+(\xi^2-1)x=0.\]
Bessel Function
\[J(1,\xi)=\sum\limits_{m=0}^{+\infty}\frac{(-1)^m}{m!(m+1)!}(\frac{\xi}{2})^{2m+1}=\frac{\xi}{2}-\frac{\xi^3}{16}+\frac{\xi^5}{384}-\frac{\xi^7}{18432}+\cdots\]
has infinitely many simple zeros in $(0,\infty)$ and satisfies
$J(1,0)=0$ when $\xi\rightarrow\,+\infty$, there holds the
asymptotical equality
\[J(1,\xi)=\sqrt{\frac{2}{\pi\xi}}\cos(\xi-\frac{3\pi}{4})+O(\xi^{-\frac{3}{2}}).\]
It implies the period function of center for system \eqref{SP},
denoted by $T(\xi,\lambda)$ which satisfies $h=G(\xi)$, it holds
\[T_{\lambda}^{\prime}(\xi,0)=0.\]
Note that $T^{\prime}(\xi)=G^{\prime}(\xi)T^{\prime}(h).$ Hence for
parameters $\xi$ and $h$ the monotonicity of the period function $T$
is the same. Thus the period function $T(\xi,\lambda)$ oscillates
and intersects the line $T(\xi,0)=2\pi$ infinitely many times, thus
it has an infinite number of roots $h_k$ such that there holds
$T(h_k,\lambda_k)=2\pi$ for all $\lambda\in\,[-1,\lambda_m)$.

The discussions mentioned as above asserts that it is possible for the existence of infinitely many critical periods for potential center \eqref{SP} as perturbation parameter $\lambda\sim\,0$. In what follows we seek to confirm it holds for the whole interval $(-1,\lambda_m).$

Without lose of generality assume that $\lambda>0$ in the following discussion, the case $\lambda\in\,[-1,0)$ can be treated similarly.

In the phase plane taking the intersection points $(x_k,0)=(2k\pi,0)$ between the oval
$\Gamma_h$ and $x-$axis, the corresponding energy level is
$h_k=G(x_k)=2(k\pi)^2$ where $k$ is a positive integer. Due to the
symmetries when $G$ is even, if $\lambda\in[0,+\infty)$ the
period of oval $\Gamma_{h_k}$ is given by
\begin{equation}\label{GJTHK}
T(h_k,\lambda)=2\sqrt{2}\int_0^{2k\pi}\frac{dx}{\sqrt{h_k-G(x)}}
              \geq\,2\sqrt{2}\int_0^{2k\pi}\frac{dx}{\sqrt{h_k-\frac{x^2}{2}}}=2\pi.
\end{equation}

Similarly consider the motion along periodic orbit passing through points $(\bar{x}_k,0)=(2k+1)\pi,0)$, the corresponding energy level is $\bar{h}_k=\frac{(2k+1)^2\pi^2}{2}+2\lambda$. The least period of the motion along the period oval is represented as
\begin{eqnarray} \nonumber \label{GJTHKB}
T(\bar{h}_k,\lambda)&=&2\sqrt{2}\int_0^{(2k+1)\pi}\frac{dx}{\sqrt{\bar{h}_k-G(x)}}\\  \nonumber
                    &=&2\sqrt{2}\int_0^{(2k+1)\pi}\frac{dx}{\sqrt{\frac{(2k+1)^2\pi^2}{2}-\frac{x^2}{2}+\lambda(1+\cos(x))}}\\ &\leq &\,2\sqrt{2}\int_0^{(2k+1)\pi}\frac{dx}{\sqrt{\frac{(2k+1)^2\pi^2}{2}-\frac{x^2}{2}}}=2\pi.
\end{eqnarray}

Combining with expressions \eqref{GJTHK} and \eqref{GJTHKB}, it can be deduced by Intermediate Value Theorem that there exists at least one energy level $h_i\in\,(h_k,\bar{h}_k)$ such that the period $T(h_i,\lambda)=2\pi$ for all $\lambda\in\,(-1,\lambda_m).$ It implies that the conservative system \eqref{SP} has countable infinitely many $2\pi$ periodic solutions. It follows from Rolle's Theorem that completes the proof of Theorem \ref{MRT}.
%\section*{Acknowledgements}


\begin{thebibliography}{99}
\bibliographystyle{amsplain}
\bibitem{GRR} G. E. Andrews, R. Askey, R. Roy; {\it Special functions}, Encyclopedia of mathematics and its applications (Book 71), Cambridge University Press; 1 edition (January 29, 2001).

\bibitem{VIAD} V. I. Arnold; {\it Loss of stability of autooscillations near a resonance and versal deformations of equivariant vector fields}. Funkts Anal. \textbf{11} (1977), 1--10.

\bibitem{JCMS} J. Chavarriga, M. Sabatini; {\it A survey of isochronous centers},
Qual. Theory Dyn. Syst. \textbf{1}(1999), 1--70.

\bibitem{CC} C. Chicone; {\it Geometric Methods for Two-Point
Nonlinear Boundary Value Problems}, J. Differential Equations \textbf{72}
(1988), 360--407.

\bibitem{ZO} Z. Opial; {\it Sur les p\'{e}riodes des solutions de l'\'{e}quation diff\'{e}rentialle $x''+g(x)=0$}, Ann. Polon. Math. \textbf{10} (1961), 49--72.

\bibitem{URABE} M. Urabe; {\it Potential forces which yield periodic motions of a fixed period}, Jour. Math. Mech. \textbf{10} (1961),
569--578.

\bibitem{SCJS} S. N. Chow, J. A. Sanders, On the number of critical points of
the period, J. Differential equations 64 (1986) 51--66.

\bibitem{CCFD} C. Chicone, F. Dumortier; {\it Finiteness for critical
periods of planar analytic vector fields}, Nonlinear Anal.
\textbf{20} (1993), 315--335.

\bibitem{MSTI} M. Sabatini, Period function's convexity for Hamiltonian centers
with separable variables, ANNALES POLONICI MATHEMATICI 85.2 (2005)
153--163.

\bibitem{CCMJ} C. Chicone, M. Jacobs; {\it Bifurcation of critical periods for plane vector fields},
Trans. Amer. Math. Soc. \textbf{312} (1989), no. 2, 433--486.

\bibitem{JF} J. P. Franc\c{o}ise; {\it The successive derivatives of the period
function of a plane vector field}, J. Differential Equations
\textbf{146} (1998), 320--335.

\bibitem{EFAA} E. Freire, A. Gasull and A. Guillamon; {\it Period function for
perturbed isochronous centres}, Qual. Theory Dyn. Syst. \textbf{3}
(2002), 275--284.

\bibitem{AJCLJY} A. Gasull, C. Liu, J. Yang; {\it On the number of critical periods for planar polynomial
systems of arbitrary degree}, J. Differential Equations \textbf{249}
(2010), no. 3, 684--692.

\bibitem{MGJV} M. Grau and J. Villadelprat; {\it Bifurcation of critical periods from Pleshkan's
isochrones}, J. London Math. Soc. \textbf{81} (2010), no. 2,
142--160.

\bibitem{DTHHZF} Ding, T., Huang, H., Zanolin, F.: A priori bounds and periodic solutions for a class of planar systems with
applications to Lotka-Volterra Equations, Disc. Cont. Dyn. Syst.
\textbf{1}, 103--117 (1995)

\bibitem{LYXZ} L. Yang, X. Zeng; {\it The period function of potential systems of polynomials with real
zeros}, Bull. Sci. Math. \textbf{133} (2009), 555--577.

\bibitem{ZZ} Z. Zhang; {\it A theorem about the differential equation $\ddot{x}+\mu
\sin(x)+x=0$  has and only has $n$ limit cycles at the domain $\mid
\dot{x} \mid \leq (n+1)\pi$}, China. Sci. \textbf{10} (1980),
941--948.

\end{thebibliography}
\end{document}